\documentclass[12pt]{article}
\usepackage[sectionbib]{natbib}
\usepackage{array,epsfig,rotating}
\usepackage[]{hyperref}
\usepackage{sectsty, secdot}
\sectionfont{\fontsize{12}{14pt plus.8pt minus .6pt}\selectfont}
\renewcommand{\theequation}{\thesection.\arabic{equation}}
\subsectionfont{\fontsize{12}{14pt plus.8pt minus .6pt}\selectfont}

\usepackage{amsmath}
\usepackage{amssymb}
\usepackage{amsfonts}
\usepackage{multirow}
\usepackage{amsthm}

\usepackage{graphicx}

\theoremstyle{definition}

\usepackage{fancyhdr}
\begin{document}


\renewcommand{\baselinestretch}{2}

\markboth{\hfill{\footnotesize\rm MIODRAG M. LOVRIC} \hfill}
{\hfill {\footnotesize\rm WHAT THE JEFFREYS--LINDLEY PARADOX REALLY IS} \hfill}

\renewcommand{\thefootnote}{}
$\ $\par


\fontsize{12}{14pt plus.8pt minus .6pt}\selectfont \vspace{0.8pc}
\centerline{\large\bf WHAT THE JEFFREYS--LINDLEY PARADOX REALLY IS:}
\vspace{2pt} 
\centerline{\large\bf CORRECTING A PERSISTENT MISCONCEPTION}
\vspace{.4cm} 
\centerline{Miodrag M. Lovric} 
\vspace{.3cm} 
\centerline{\it Department of Mathematics and Statistics}
\centerline{\it Radford University}
\centerline{\it P.O. Box 6942, Radford, VA 24142, USA}
\centerline{\it E-mail: mlovric@radford.edu}
 \vspace{.55cm} \fontsize{9}{11.5pt plus.8pt minus.6pt}\selectfont


\begin{quotation}
\noindent {\it Abstract:}
The Jeffreys--Lindley paradox stands as the most profound divergence between frequentist and Bayesian approaches to hypothesis testing. Yet despite more than six decades of discussion, this paradox remains frequently misunderstood---even in the pages of leading statistical journals. In a 1993 paper published in \emph{Statistica Sinica}, \citet{Robert1993} characterized the Jeffreys--Lindley paradox as ``the fact that a point null hypothesis will always be accepted when the variance of a conjugate prior goes to infinity.'' This characterization, however, describes a different phenomenon entirely---what we term Bartlett's Anomaly---rather than the Jeffreys--Lindley paradox as originally formulated. The paradox, as presented by \citet{Lindley1957}, concerns what happens as sample size increases without bound while holding the significance level fixed, not what happens as prior variance diverges. This distinction is not merely terminological: the two phenomena have different mathematical structures, different implications, and require different solutions. The present paper aims to clarify this confusion, demonstrating through Lindley's own equations that he was concerned exclusively with sample size asymptotics. We show that even Jeffreys himself underestimated the practical frequency of the paradox. Finally, we argue that the only genuine resolution lies in abandoning point null hypotheses in favor of interval nulls, a paradigm shift that eliminates the paradox and restores harmony between Bayesian and frequentist inference.

\vspace{9pt}
\par
\end{quotation}\par

\def\thefigure{\arabic{figure}}
\def\thetable{\arabic{table}}

\renewcommand{\theequation}{\thesection.\arabic{equation}}

\fontsize{12}{14pt plus.8pt minus .6pt}\selectfont

\section{Introduction}
\setcounter{equation}{0}

The Jeffreys--Lindley paradox occupies a unique position in the foundations of statistical inference. It is, as \citet{Lovric2019} observed, the ``Gordian Knot'' of modern statistics---a fundamental tension that divides frequentist and Bayesian camps in seemingly irreconcilable ways. For more than sixty years, statisticians have debated its meaning, its implications, and its resolution, yet no consensus has emerged \citep[for a comprehensive overview, see][]{Lovric2025a}. Indeed, \citet{GelmanShalizi2013} concluded that ``the Jeffreys--Lindley paradox... is really a problem without a solution.''

Given the paradox's complexity, it is perhaps unsurprising that misunderstandings have proliferated. What is surprising---and troubling---is that some of these misunderstandings have appeared in influential papers by distinguished statisticians. The present paper addresses one such case: the characterization of the Jeffreys--Lindley paradox offered by \citet{Robert1993} in \emph{Statistica Sinica}.

Robert, author of the celebrated textbook \emph{The Bayesian Choice} \citep{Robert2007}, is unquestionably one of the leading figures in modern Bayesian statistics. Yet in his 1993 paper ``A Note on Jeffreys--Lindley Paradox,'' he opens with the following definition:
\begin{quotation}
``The Jeffreys--Lindley paradox, namely the fact that a point null hypothesis will always be accepted when the variance of a conjugate prior goes to infinity, has often been argued to imply prohibiting the use of improper priors in hypothesis testing.'' \citep[p.~601]{Robert1993}
\end{quotation}
This characterization, we shall demonstrate, describes a phenomenon entirely distinct from what Lindley analyzed in his seminal 1957 paper. Robert's ``paradox'' concerns the behavior of Bayes factors as prior variance $\sigma_0^2 \to \infty$ with \emph{fixed} data. The Jeffreys--Lindley paradox as Lindley originally formulated it concerns the behavior of posterior probabilities as sample size $n \to \infty$ with a \emph{fixed} significance level. These are fundamentally different asymptotic regimes with different mathematical structures and different practical implications.

The distinction matters for several reasons. First, it affects how we diagnose the source of the paradox and therefore how we might resolve it. Second, the two phenomena have been conflated in the subsequent literature, propagating confusion---a point also noted in the comprehensive historical overview by \citet{WagenmakersLy2022}. Third, solving one problem does not solve the other---as Robert himself later acknowledged when he abandoned his 1993 ``solution'' as ``flawed from the measure-theoretic angle'' \citep[p.~137]{RousseauRobert2011}.

The goals of this paper are threefold: (1) to clarify definitively what the Jeffreys--Lindley paradox actually is, by returning to Lindley's original formulation; (2) to distinguish it sharply from Bartlett's Anomaly, which Robert addressed; and (3) to explain why the only genuine resolution lies in a paradigm shift from point null hypotheses to interval nulls. Along the way, we note that even Harold Jeffreys, who first identified the tension between $p$-values and Bayes factors, appears to have underestimated its practical significance.

\section{Lindley's Original Formulation}
\setcounter{equation}{0}

To understand what the Jeffreys--Lindley paradox truly is, we must return to its source. \citet{Lindley1957} considered testing a point null hypothesis $H_0: \theta = \theta_0$ within a normal model $X \sim N(\theta, \sigma^2)$ with known variance. Following the approach pioneered by \citet{Jeffreys1939}, Lindley employed a mixed prior distribution that assigns probability mass $c$ to the null value $\theta_0$ and distributes the remainder $(1-c)$ uniformly over an interval of width $I$ containing $\theta_0$. (Note that the variance of this uniform distribution is $I^2/12$, so increasing $I$ is equivalent to increasing the prior variance under $H_1$.) The choice of uniform prior is not essential to the paradox; any proper prior assigning positive mass to the point null will produce the same qualitative behavior. This ``spike-and-slab'' prior has the form:
\begin{equation}
P(\theta) = c\,\delta_{\theta=\theta_0} + (1-c)\,g(\theta)\,\mathbf{1}_{\theta \neq \theta_0},
\end{equation}
where $\delta_{\theta=\theta_0}$ denotes the Dirac mass at $\theta_0$.

Lindley then supposed that the sample mean $\bar{x}$ was ``just significant'' at the $\alpha$ level, meaning
\begin{equation}
\bar{x} = \theta_0 + \lambda_{\alpha/2} \frac{\sigma}{\sqrt{n}},
\end{equation}
where $\lambda_{\alpha/2}$ is the upper $\alpha/2$ quantile of the standard normal distribution. Under these conditions, Lindley derived the posterior probability of the null hypothesis as:
\begin{equation}\label{eq:Lindley}
P(H_0 \mid \text{just significant } \bar{x}) = \left[ c\exp(-\lambda_{\alpha/2}^2/2) + (1-c)\sigma\sqrt{\frac{2\pi}{n}} \right]^{-1} \cdot c\exp(-\lambda_{\alpha/2}^2/2).
\end{equation}

\textbf{The crucial observation is this: Equation (\ref{eq:Lindley}) contains the sample size $n$, the significance level $\alpha$ (through $\lambda_{\alpha/2}$), and the prior probability $c$---but it contains no term for the prior variance under the alternative hypothesis.} Lindley was not analyzing what happens when prior variance increases; he was analyzing what happens when sample size increases while the $p$-value remains fixed at $\alpha$. In plain terms: Lindley's paradox is about \emph{more data}, not \emph{vaguer priors}.

From Equation (\ref{eq:Lindley}), Lindley deduced that as $n \to \infty$, the posterior probability $P(H_0 \mid \bar{x}) \to 1$. This occurs because the term $(1-c)\sigma\sqrt{2\pi/n}$ vanishes, leaving only the contribution from the point mass at $\theta_0$. Thus, for any fixed significance level $\alpha$ and any fixed prior probability $c > 0$, there exists a sample size $n$ large enough that:
\begin{enumerate}
\item The sample mean $\bar{x}$ is statistically significant at the $\alpha$ level (frequentist conclusion: reject $H_0$), and simultaneously
\item The posterior probability that $\theta = \theta_0$ exceeds $(1-\alpha)$ (Bayesian conclusion: strong support for $H_0$).
\end{enumerate}
Lindley called this conflicting situation ``the paradox'' (p.~187) and ``strong contrast'' (p.~190).

Although the derivation above uses the normal model for concreteness, the distinction between sample-size asymptotics and prior-diffuseness holds much more generally. The paradox arises whenever: (i) the likelihood concentrates at rate $n^{-1/2}$ (or more generally, the posterior variance shrinks with $n$), and (ii) a point null receives positive prior mass. These conditions are satisfied by most regular parametric models; see \citet{BergerDelampady1987} for discussion of non-normal cases.

\section{Bartlett's Correction and a Different Phenomenon}
\setcounter{equation}{0}

Shortly after Lindley's paper appeared, \citet{Bartlett1957} published a brief comment correcting a slip in Lindley's analysis. Bartlett noted that the uniform density over the interval $I$ should contribute a factor of $1/I$, giving the corrected posterior probability:
\begin{equation}\label{eq:Bartlett}
\bar{c} = \left[ c\exp(-\lambda_{\alpha/2}^2/2) + \frac{(1-c)}{I}\sigma\sqrt{\frac{2\pi}{n}} \right]^{-1} \cdot c\exp(-\lambda_{\alpha/2}^2/2).
\end{equation}

Bartlett observed that this correction reveals a dependence on the interval width $I$ (equivalently, the prior variance under the alternative). He noted wryly that ``one might be tempted to put $I$ infinity [but] the silly answer $\bar{c} = 1$ ensues'' (p.~533). In plain terms: Bartlett's Anomaly is about \emph{vaguer priors}, not \emph{more data}.

This observation exposes a troubling fact: for any fixed data (and hence fixed $\lambda_{\alpha/2}$), the posterior probability $P(H_0 \mid \bar{x}) \to 1$ as the prior becomes more diffuse ($I \to \infty$ or equivalently $\sigma_0^2 \to \infty$). This means that a data analyst can ``prove'' any null hypothesis simply by choosing a sufficiently vague prior under the alternative---a form of what we might call ``posterior probability hacking.''

\textbf{This is a genuine problem, but it is not the Jeffreys--Lindley paradox.} It is a different phenomenon with a different driver:

\begin{center}
\begin{tabular}{lcc}
\hline
& \textbf{Jeffreys--Lindley Paradox} & \textbf{Bartlett's Anomaly} \\
\hline
What varies & Sample size $n \to \infty$ & Prior variance $\sigma_0^2 \to \infty$ \\
What is fixed & Prior, significance level $\alpha$ & Data, sample size $n$ \\
Driver & Likelihood concentration & Prior diffuseness \\
Result & $P(H_0 \mid x) \to 1$ & $P(H_0 \mid x) \to 1$ \\
\hline
\end{tabular}
\end{center}

\begin{figure}[htbp]
\centering
\includegraphics[width=\textwidth]{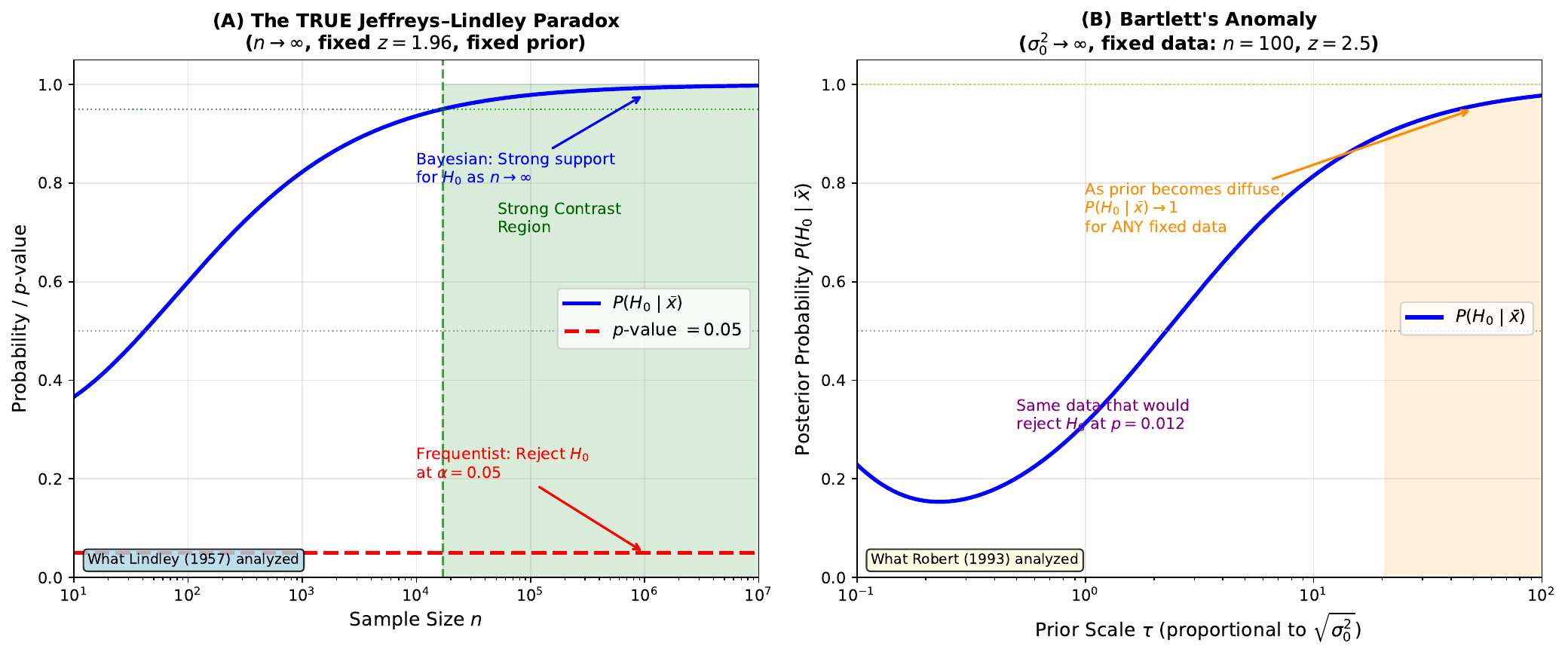}
\caption{Visual comparison of two distinct phenomena, computed for a normal model with known variance. \textbf{(A)} The Jeffreys--Lindley paradox (Lindley, 1957) occurs when sample size $n \to \infty$ with fixed prior and fixed $p$-value. Here $z = 1.96$ (corresponding to $\alpha = 0.05$), prior scale $\tau = \sigma_0/\sigma = 1$ (a standard ``unit information'' choice), and prior probability $\pi_0 = 0.5$ (equal prior odds). Note that $\bar{x}$ approaches $\theta_0$ as $n \to \infty$ to maintain the fixed-$z$ constraint; the observed effect shrinks but remains ``significant.'' The posterior probability $P(H_0|\bar{x})$ increases toward 1 as $n$ grows, creating ``strong contrast'' with the frequentist rejection. \textbf{(B)} Bartlett's Anomaly (Bartlett, 1957) occurs when prior variance increases with fixed data and fixed sample size. Here $n = 100$, $z = 2.5$ (yielding $p \approx 0.012$), and $\pi_0 = 0.5$; the data $\bar{x}$ remain unchanged throughout. The horizontal axis shows prior scale $\tau = \sigma_0/\sigma$ (so prior variance is $\tau^2\sigma^2$). Even data yielding $p = 0.012$ leads to $P(H_0|\bar{x}) \to 1$ as the prior becomes diffuse. These are fundamentally different phenomena with different drivers: likelihood concentration (A) versus prior diffuseness (B).}
\label{fig:comparison}
\end{figure}

While both phenomena lead to the posterior favoring $H_0$, they arise from entirely different mechanisms. The Jeffreys--Lindley paradox is driven by the geometry of the likelihood function as it concentrates around the null value at rate $1/\sqrt{n}$. Bartlett's Anomaly is driven by the spreading of prior mass under the alternative, which penalizes $H_1$ regardless of what the data show.

Why has this confusion persisted for so long? We suggest that the conflation is cognitively natural: both phenomena share the same symptom ($P(H_0 \mid x) \to 1$), and both involve a parameter becoming ``large'' or ``diffuse''---sample size in one case, prior variance in the other. Without careful attention to \emph{what is varying} and \emph{what is held fixed}, the two are easily mistaken for one another.

We propose calling this latter phenomenon \textbf{Bartlett's Anomaly} rather than ``Bartlett's paradox'' because there is nothing truly paradoxical about it: diffuse priors on continuous parameters lead to improper Bayes factors, and improper procedures yield nonsensical answers. This has been well understood since at least \citet{DeGroot1982}.

\section{Reinterpreting Robert (1993)}
\setcounter{equation}{0}

With this distinction in hand, we can now reexamine Robert's 1993 paper. His opening sentence defines the Jeffreys--Lindley paradox as ``the fact that a point null hypothesis will always be accepted when the variance of a conjugate prior goes to infinity.'' This description corresponds to Bartlett's Anomaly rather than the Jeffreys--Lindley paradox as Lindley originally formulated it.

The key observation supporting this interpretation is that \textbf{Lindley's Equation (\ref{eq:Lindley}) does not contain any term for prior variance}. Lindley was not even aware of the dependence on prior variance until Bartlett's correction---as evidenced by his omission of the factor $1/I$ from his original formula. The paradox Lindley described was entirely about sample size asymptotics.

Robert's paper then proceeds to develop an ingenious ``solution'' in which the prior probability $\rho_0$ of the null hypothesis is made to depend on the prior variance $\sigma^2$ under the alternative, specifically requiring that:
\begin{equation}
\frac{1 - \rho_0(\sigma)}{\sigma} = \text{constant}.
\end{equation}
This ensures that as $\sigma \to \infty$, the prior probability of $H_0$ also approaches appropriate limits, yielding a ``noninformative answer'' that is ``surprisingly close to the classical $p$-value'' (p.~602).

This is mathematically elegant, but it addresses Bartlett's Anomaly rather than the Jeffreys--Lindley paradox. It resolves the prior-diffuseness phenomenon by preventing the prior from becoming arbitrarily diffuse without compensating adjustment to the prior odds. It does \emph{not} address the sample-size phenomenon, which occurs even with \emph{fixed}, proper priors as sample size increases.

Three observations support this interpretation:
\begin{enumerate}
\item If one reads ``Jeffreys--Lindley paradox'' as ``Bartlett's Anomaly'' throughout Robert's paper, the analysis is entirely coherent and correct.

\item Robert himself later moved away from his 1993 approach, acknowledging in \citet[p.~137]{RousseauRobert2011} that it was ``flawed from the measure-theoretic angle.''

\item The subsequent literature has sometimes followed Robert's usage. For example, \citet[p.~12290]{VillaWalker2017} describe Lindley as showing that ``if the prior for the location parameter, in the alternative model to the parameter being zero, has infinite variance, then the Bayesian will always select the null model, regardless of the observed data.'' However, Lindley's original analysis held prior variance fixed and varied sample size.
\end{enumerate}

We emphasize that Robert's 1993 paper made a genuine contribution: it provided a coherent treatment of how prior odds should scale with increasing prior variance to avoid the pathologies of diffuse priors. The issue is terminological rather than mathematical---the phenomenon Robert analyzed is real and important, but it is distinct from what Lindley originally called ``the paradox.''

\section{The Practical Dimension: Jeffreys Also Underestimated}
\setcounter{equation}{0}

Having clarified what the Jeffreys--Lindley paradox is, we now address its practical significance. Harold Jeffreys was the first to notice the potential disagreement between $p$-values and Bayes factors. In the Appendix to his \emph{Theory of Probability}, he wrote:
\begin{quotation}
``At large numbers of observations there is a difference, since the tests based on the integral [$p$-value] would sometimes assert significance at departures that would actually give $K > 1$ [Bayes factor supports $H_0$]. Thus there may be opposite decisions in such cases. \emph{But these will be very rare}.'' \citep[p.~435, emphasis added]{Jeffreys1961}
\end{quotation}

Jeffreys believed that while the paradox existed in principle, it would ``extremely rarely'' affect practical decisions. This assessment appears to have influenced the broader perception of the paradox as a theoretical curiosity of little applied relevance---a view echoed by \citet[p.~322]{BergerDelampady1987}, who suggested the paradox might be of ``questionable relevance.''

However, simulation studies in \citet{Lovric2019} demonstrate that Jeffreys' assessment was overly optimistic. The paradox is not rare; it occurs regularly whenever sample sizes are sufficiently large. Table~\ref{tab:frequency} shows the minimum sample size required for the paradox to manifest in its ``strong contrast'' form---Lindley's term for the situation where the $p$-value equals $\alpha$ (frequentist rejects $H_0$) while simultaneously the posterior probability $P(H_0 \mid \bar{x}) > 1-\alpha$ (Bayesian strongly supports $H_0$). These values were computed by solving Equation~(\ref{eq:Lindley}) for $n$ given the constraint $P(H_0 \mid \bar{x}) = 1 - \alpha$, with prior probability $c = 0.5$ for Lindley's setup. For the normal conjugate case, the Bayes factor formula $B_{01} = (1 + n\tau^2)^{-1/2} \exp(z^2 n\tau^2 / [2(1 + n\tau^2)])$ was used with $\tau = 1$ and $z = \lambda_{\alpha/2}$. The choice $\tau = 1$ (prior standard deviation equal to sampling standard deviation) represents a ``unit information prior'' \citep{Kass1995}---a principled default that assigns one observation's worth of information to the prior, widely used in objective Bayesian analysis.

\begin{table}[h]
\caption{Minimum sample size for the Jeffreys--Lindley paradox to occur in strong contrast form.}\label{tab:frequency}
\centering
\begin{tabular}{ccrr}
\hline
$\alpha$ level & $P(H_0 \mid \text{just significant } \bar{x})$ & Lindley's setup & Normal conjugate \\
\hline
0.050 & 0.950 & 105,685 & 16,816 \\
0.040 & 0.960 & 245,701 & 39,098 \\
0.030 & 0.970 & 728,954 & 116,011 \\
0.020 & 0.980 & 3,380,074 & 537,945 \\
0.010 & 0.990 & 46,875,786 & 2,195,961 \\
0.005 & 0.995 & 657,481,111 & 104,625,626 \\
\hline
\end{tabular}
\end{table}

In the Zettabyte Era, datasets routinely contain millions or billions of observations. Genomics studies analyze millions of genetic variants; social media platforms process billions of user interactions; large-scale clinical trials and electronic health records encompass populations that dwarf traditional sample sizes. At such scales, the paradox is not a theoretical curiosity but an everyday occurrence. The claim by \citet[p.~235]{EdwardsLindmanSavage1963} that ``results of [Bayesian and classical testing] procedures will usually agree'' is simply untenable for modern large-scale data analysis.

\section{Toward Resolution: Interval Null Hypotheses}
\setcounter{equation}{0}

If the Jeffreys--Lindley paradox cannot be resolved by adjusting priors (as Robert attempted), what hope remains? The answer lies not in modifying the inference procedure but in reformulating the hypothesis itself.

The root of the paradox is the assignment of positive prior probability to a single point $\theta_0$ in a continuous parameter space. This point has Lebesgue measure zero, yet it receives finite prior mass---a measure-theoretic peculiarity that Jeffreys himself recognized was necessary to make Bayesian hypothesis testing ``competitive'' with frequentist methods. As the likelihood concentrates around $\bar{x}$ at rate $1/\sqrt{n}$, and $\bar{x}$ is constrained to approach $\theta_0$ by the fixed-$p$ condition, the point mass at $\theta_0$ captures an ever-increasing share of the posterior.

The resolution is to replace the point null $H_0: \theta = \theta_0$ with an \textbf{interval null} (or ``region of practical equivalence''):
\begin{equation}
H_0: |\theta - \theta_0| \leq \delta \quad \text{vs.} \quad H_1: |\theta - \theta_0| > \delta,
\end{equation}
where $\delta > 0$ represents the smallest scientifically meaningful effect size.

Under this formulation, both $H_0$ and $H_1$ are composite hypotheses with positive Lebesgue measure. Using continuous priors $\pi_0$ on $[\theta_0 - \delta, \theta_0 + \delta]$ and $\pi_1$ on the complement, the Bayes factor becomes:
\begin{equation}
B_{01} = \frac{\int_{|\theta - \theta_0| \leq \delta} L(\theta) \pi_0(\theta) \, d\theta}{\int_{|\theta - \theta_0| > \delta} L(\theta) \pi_1(\theta) \, d\theta}.
\end{equation}

Under the ``just significant'' condition where $\bar{x} = \theta_0 + z\sigma/\sqrt{n}$ with fixed $z$, the observed effect shrinks toward zero as $n$ increases, eventually falling well within the equivalence region $[\theta_0 - \delta, \theta_0 + \delta]$. In this regime, the Bayes factor $B_{01} \to \infty$, strongly supporting $H_0$. But crucially, this no longer conflicts with frequentist inference: under equivalence testing \citep{Schuirmann1987}, the frequentist also concludes in favor of $H_0$ when the observed effect is practically negligible.

The conflict disappears not because the Bayes factor behaves differently, but because the frequentist framework changes. Under point-null testing at fixed $\alpha$, the frequentist rejects $H_0$ whenever $|z| > z_{\alpha/2}$, regardless of sample size---this creates the paradox. Under equivalence testing, the frequentist asks whether the data support $|\theta - \theta_0| \leq \delta$, and when the observed effect falls within this region, the equivalence test concludes in favor of $H_0$. Both paradigms now agree.

While the interval null formulation dissolves the theoretical paradox, it introduces the practical task of specifying $\delta$---the region of practical equivalence---which must be justified by subject-matter expertise. Crucially, however, $\delta$ has direct scientific meaning (the smallest effect of practical importance), unlike the arbitrary prior variance that drives Bartlett's Anomaly.

We acknowledge that other Bayesian approaches have been proposed to address issues related to hypothesis testing with point nulls. These include nonlocal priors \citep{JohnsonRossell2010}, which place zero density at the null value under $H_1$; intrinsic Bayes factors \citep{BergerPericchi1996}, which use training samples to construct ``automatic'' reference priors; and fractional Bayes factors, which use a fraction of the likelihood for prior calibration. Each of these methods has merit and effectively addresses Bartlett's Anomaly by preventing arbitrary posterior accumulation from diffuse priors. However, they do not resolve the Jeffreys--Lindley paradox proper, which persists even under fixed, proper priors as $n \to \infty$. We contend that interval nulls offer a more fundamental resolution because they address the core issue: point nulls ask a question (``is $\theta$ exactly equal to $\theta_0$?'') that no finite dataset can meaningfully answer. Interval nulls reformulate the question in scientifically meaningful terms.

As \citet{Lovric2025b} observed, this represents a genuine dissolution of the paradox: ``the only solution of JL paradox is its dissolution.'' The paradox does not occur when inference is conducted on interval null hypotheses because both paradigms are then answering the same question: whether the data are consistent with a practically negligible effect.

\section{Concluding Remarks}
\setcounter{equation}{0}

The Jeffreys--Lindley paradox has puzzled statisticians for over sixty years. Much of this puzzlement stems from confusion about what the paradox actually is. The present paper has sought to clarify this definitively:

\begin{enumerate}
\item \textbf{The Jeffreys--Lindley paradox} as originally formulated concerns the behavior of posterior probabilities as sample size $n \to \infty$ with a fixed significance level. Lindley's original equation contains no prior variance term; his analysis was entirely about sample size asymptotics.

\item \textbf{Bartlett's Anomaly} concerns the behavior of Bayes factors as prior variance $\sigma_0^2 \to \infty$ with fixed data. This is a different phenomenon with different mathematical structure.

\item \textbf{Robert (1993)} addressed Bartlett's Anomaly---a genuine and important phenomenon---while using the terminology ``Jeffreys--Lindley paradox.'' We suggest these two phenomena be clearly distinguished in future work.

\item \textbf{Even Jeffreys underestimated} the practical frequency of the paradox, claiming conflicting decisions would be ``very rare.'' Modern simulation shows they are common with large samples.

\item \textbf{A natural resolution} is a paradigm shift from point null hypotheses to interval nulls. Under this formulation, the paradox dissolves because both Bayesian and frequentist frameworks answer the same scientifically meaningful question---whether the effect is practically negligible---and reach the same conclusion.
\end{enumerate}

We hope this clarification will benefit readers of \emph{Statistica Sinica} and the broader statistical community. The Jeffreys--Lindley paradox need not divide Bayesians and frequentists; properly understood, it points toward a deeper unity. When hypotheses encode \emph{practical} rather than exact equality, all coherent frameworks converge. The conflict vanishes not by altering probability calculus, but by refining the questions we ask.

\section*{Acknowledgements}

The author thanks colleagues for helpful discussions. This paper is offered in a spirit of scholarly clarification, with deep respect for the contributions of Christian Robert and all researchers who have grappled with this fundamental problem.
\par


\bibhang=1.7pc
\bibsep=2pt
\fontsize{9}{14pt plus.8pt minus .6pt}\selectfont
\renewcommand\bibname{\large \bf References}
\expandafter\ifx\csname
natexlab\endcsname\relax\def\natexlab#1{#1}\fi
\expandafter\ifx\csname url\endcsname\relax
  \def\url#1{\texttt{#1}}\fi
\expandafter\ifx\csname urlprefix\endcsname\relax\def\urlprefix{URL}\fi

\vskip .65cm
\noindent
Department of Mathematics and Statistics, Radford University, Radford, VA 24142, USA
\vskip 2pt
\noindent
E-mail: mlovric@radford.edu

\newpage
\appendix
\section{Technical Derivations}\label{app:derivations}
\setcounter{equation}{0}
\renewcommand{\theequation}{A.\arabic{equation}}

This appendix provides detailed derivations of the key equations discussed in the main text.

\subsection{Derivation of Lindley's Formula}

Consider testing $H_0: \theta = \theta_0$ versus $H_1: \theta \neq \theta_0$ based on a sample mean $\bar{x}$ from $N(\theta, \sigma^2/n)$ with known $\sigma^2$. The mixed prior assigns probability $c$ to $H_0$ and spreads probability $(1-c)$ uniformly over an interval $[\theta_0 - I/2, \theta_0 + I/2]$ of fixed width $I$.

\textbf{Step 1: The correct derivation with proper normalization.}

The posterior odds are:
\begin{equation}
\frac{P(H_0 \mid \bar{x})}{P(H_1 \mid \bar{x})} = \frac{c}{1-c} \cdot \frac{f(\bar{x} \mid \theta_0)}{\int_{-I/2}^{I/2} f(\bar{x} \mid \theta_0 + u) \cdot \frac{1}{I} \, du},
\end{equation}
where $f(\bar{x} \mid \theta) = \sqrt{n/(2\pi\sigma^2)} \exp(-n(\bar{x}-\theta)^2/(2\sigma^2))$.

Suppose $\bar{x}$ is ``just significant'' at level $\alpha$, so $\bar{x} = \theta_0 + \lambda_{\alpha/2}\sigma/\sqrt{n}$ where $\lambda_{\alpha/2} = \Phi^{-1}(1-\alpha/2)$. Then:
\begin{equation}
f(\bar{x} \mid \theta_0) = \sqrt{\frac{n}{2\pi\sigma^2}} \exp\left(-\frac{\lambda_{\alpha/2}^2}{2}\right).
\end{equation}

For the integral in the denominator, as $n$ becomes large, the likelihood concentrates near $\bar{x}$. Provided $I$ is fixed and large enough to contain the region of substantial likelihood, Laplace approximation gives:
\begin{equation}
\int_{-I/2}^{I/2} f(\bar{x} \mid \theta_0 + u) \cdot \frac{1}{I} \, du \approx \frac{1}{I} \cdot \sqrt{\frac{2\pi\sigma^2}{n}} \cdot f(\bar{x} \mid \bar{x}) = \frac{1}{I} \cdot \sqrt{\frac{2\pi\sigma^2}{n}} \cdot \sqrt{\frac{n}{2\pi\sigma^2}} = \frac{1}{I}.
\end{equation}

Thus, with proper normalization, the marginal likelihood under $H_1$ is approximately $1/I$, and the posterior odds are:
\begin{equation}\label{eq:correct}
\frac{P(H_0 \mid \bar{x})}{P(H_1 \mid \bar{x})} \approx \frac{c}{1-c} \cdot I \cdot f(\bar{x} \mid \theta_0) = \frac{c}{1-c} \cdot I \cdot \sqrt{\frac{n}{2\pi\sigma^2}} \exp\left(-\frac{\lambda_{\alpha/2}^2}{2}\right).
\end{equation}

This expression depends on both $n$ and $I$.

\textbf{Step 2: Why Lindley's published formula differs.}

In Lindley's original presentation, the uniform density $1/I$ was treated as an unspecified constant $k$ rather than being explicitly tracked. This was likely intentional: Lindley was focused on the behavior as $n \to \infty$ with $I$ fixed and ``sufficiently large,'' so the exact value of $I$ seemed immaterial to his analysis. However, this choice had the unintended consequence of obscuring the formula's dependence on prior variance, which later led to the conflation we address in this paper.

Working through the algebra with this implicit constant leads to a term proportional to $\sigma\sqrt{2\pi/n}$ appearing in the denominator rather than $1/I$. The result is:
\begin{equation}\label{eq:Lindley_appendix}
P(H_0 \mid \bar{x}) = \frac{c \exp(-\lambda_{\alpha/2}^2/2)}{c \exp(-\lambda_{\alpha/2}^2/2) + (1-c)\sigma\sqrt{2\pi/n}},
\end{equation}
which is Equation~(2.3) in the main text. \textbf{This expression depends on $n$ but not on $I$}, because the factor $1/I$ was absorbed into an implicit constant.

As $n \to \infty$, the term $(1-c)\sigma\sqrt{2\pi/n} \to 0$, so $P(H_0 \mid \bar{x}) \to 1$---this is the Jeffreys--Lindley paradox. The key point for our purposes is that Lindley's analysis concerns the $n \to \infty$ limit with $I$ fixed (and implicitly absorbed), not the $I \to \infty$ limit.

\subsection{Bartlett's Correction}

\citet{Bartlett1957} noted that Lindley's formula should explicitly include the factor $1/I$ from the uniform prior density. The corrected posterior probability is:
\begin{equation}
P(H_0 \mid \bar{x}) = \frac{c \exp(-\lambda_{\alpha/2}^2/2)}{c \exp(-\lambda_{\alpha/2}^2/2) + \frac{(1-c)}{I}\sigma\sqrt{2\pi/n}}.
\end{equation}

Now if we hold $n$ and the data fixed but let $I \to \infty$, the second term in the denominator vanishes, giving $P(H_0 \mid \bar{x}) \to 1$. This is Bartlett's Anomaly: diffuse priors favor $H_0$ regardless of the data. Bartlett's correction reveals the dependence on $I$ that Lindley's original formula obscured.

\subsection{Interval Null Asymptotics: Why the Paradox Dissolves}

For interval nulls $H_0: |\theta - \theta_0| \leq \delta$ versus $H_1: |\theta - \theta_0| > \delta$ with continuous priors $\pi_0$ and $\pi_1$, the Bayes factor is:
\begin{equation}
B_{01} = \frac{\int_{\theta_0-\delta}^{\theta_0+\delta} L_n(\theta) \pi_0(\theta) \, d\theta}{\int_{|\theta-\theta_0|>\delta} L_n(\theta) \pi_1(\theta) \, d\theta},
\end{equation}
where $L_n(\theta) \propto \exp(-n(\bar{x}-\theta)^2/(2\sigma^2))$.

Under the ``just significant'' condition with $\bar{x} = \theta_0 + z\sigma/\sqrt{n}$ for fixed $z$, we have $\bar{x} \to \theta_0$ as $n \to \infty$. For sufficiently large $n$, the maximum likelihood estimate $\hat{\theta} = \bar{x}$ lies inside the interval null region $[\theta_0 - \delta, \theta_0 + \delta]$.

The key observation is that the numerator and denominator scale very differently:

\textbf{Numerator (integral over $|\theta - \theta_0| \leq \delta$):} The likelihood peaks at $\hat{\theta} = \bar{x}$, which lies inside this region. By Laplace approximation:
\begin{equation}
\int_{\theta_0-\delta}^{\theta_0+\delta} L_n(\theta) \pi_0(\theta) \, d\theta \approx L_n(\bar{x}) \cdot \pi_0(\bar{x}) \cdot \sqrt{\frac{2\pi\sigma^2}{n}} = O(n^{-1/2}).
\end{equation}

\textbf{Denominator (integral over $|\theta - \theta_0| > \delta$):} The likelihood must be evaluated at points at least distance $\delta$ from $\theta_0$, and hence approximately distance $\delta$ from $\bar{x}$ for large $n$. At such points:
\begin{equation}
L_n(\theta_0 \pm \delta) \propto \exp\left(-\frac{n(\bar{x} - \theta_0 \mp \delta)^2}{2\sigma^2}\right) = \exp\left(-\frac{n\delta^2}{2\sigma^2} + O(\sqrt{n})\right).
\end{equation}
This decays \emph{exponentially} in $n$, not polynomially.

Therefore, the Bayes factor behaves as:
\begin{equation}
B_{01} \sim \frac{O(n^{-1/2})}{\exp(-n\delta^2/(2\sigma^2))} \to \infty \quad \text{as } n \to \infty.
\end{equation}

\textbf{Why this dissolves the paradox:} The Jeffreys--Lindley paradox arose from a \emph{conflict} between frequentist and Bayesian conclusions: the frequentist rejected $H_0$ (based on a significant $p$-value) while the Bayesian strongly supported $H_0$ (based on $B_{01} \to \infty$). With interval nulls, this conflict disappears because the frequentist framework changes:

\begin{itemize}
\item Under \textbf{point-null testing at fixed $\alpha$}: The frequentist rejects $H_0$ whenever $|z| > z_{\alpha/2}$, regardless of $n$. This creates the paradox.

\item Under \textbf{equivalence testing} \citep{Schuirmann1987}: The frequentist asks whether the data support $|\theta - \theta_0| \leq \delta$. When the observed effect $\bar{x} - \theta_0 = z\sigma/\sqrt{n} \to 0$ falls well within the equivalence region, the equivalence test concludes in favor of $H_0$.
\end{itemize}

Thus, both Bayesian and frequentist procedures agree: when the observed effect is practically negligible (inside the equivalence region), both provide strong support for $H_0$. The paradox is dissolved not because the Bayes factor has a finite limit, but because \emph{both frameworks now answer the same scientifically meaningful question} and reach the same conclusion.

\end{document}